\documentclass{amsart}

\usepackage{amssymb,color}
\usepackage{amsfonts}
\usepackage{amsmath}
\usepackage{euscript}
\usepackage{enumerate}
\usepackage{graphics}

\usepackage{hyperref}
\usepackage{pdfsync}
\newtheorem*{Theorem1'}{Theorem 1'}

\theoremstyle{definition}

\theoremstyle{remark}

\newcommand \Z{{\mathbb Z}}

\newcommand \N{{\mathbb N}}

\def\d{\sqrt{d}}

\begin{document}

\title[Closed formulae for certain Fermat-Pell equations]{Closed formulae for certain Fermat-Pell equations}

\author{Fernando Szechtman}
\address{Department of Mathematics and Statistics, University of Regina, Canada}
\email{fernando.szechtman@gmail.com}
\thanks{This research was partially supported by an NSERC grant}

\subjclass[2010]{11Y65, 11A55, 11D09, 11Y50, 11B37}



\keywords{Fermat-Pell equation, continued fraction}

\begin{abstract} Given positive integers $j,k$, with $j\geq 2$, we show that there are positive integers $d,e$ such that $\d$
has continued fraction expansion $\d=[e,\overline{k,\dots,k,2e}]$, with period $j$, if and only if $k$ is even or $3\nmid j$,
in which case we give closed formulae to find all such $d,e$ as well as the smallest solution in positive integers to the Fermat-Pell
equation $X^2-dY^2=(-1)^j$.
\end{abstract}

\maketitle


\section{Introduction}

Given a positive integer $d$ that is not a perfect square, its square root $\sqrt{d}$ has an infinite periodic continued fraction expansion
$$
\d=[e,\overline{a_1,\dots,a_{j-1},2e}],
$$
where $e=[\d]$ is the integral part of $\d$, $j\geq 1$ is the period of this expansion, and
$$
a_1=a_{j-1},\; a_2=a_{j-2},\;\dots
$$
are positive integers \cite{S}. Set
\begin{equation}
\label{zero}
\left(
  \begin{array}{cc}
    y & x \\
    z & w \\
  \end{array}
\right)=\left(
  \begin{array}{cc}
    0 & 1 \\
    1 & 2e \\
  \end{array}
\right)\left(\begin{array}{cc}
    0 & 1 \\
    1 & a_{j-1} \\
  \end{array}
\right)\cdots
\left(\begin{array}{cc}
    0 & 1 \\
    1 & a_1 \\
  \end{array}
\right)\left(\begin{array}{cc}
    0 & 1 \\
    1 & e \\
  \end{array}
\right).
\end{equation}
Then $(x,y)\in\N\times \N$ is a solution to the Fermat-Pell equation
\begin{equation}
\label{uno}
X^2-dY^2=(-1)^j.
\end{equation}
Moreover, if $j$ is even, then
\begin{equation}
\label{dos}
X^2-dY^2=-1
\end{equation}
has no integral solutions, and $(x_1,y_1)\in\Z\times\Z$ is a solution to
\begin{equation}
\label{tres}
X^2-dY^2=1
\end{equation}
if and only if
\begin{equation}
\label{cuatro}
x_1+y_1\d=\pm (x+y\d)^k,\quad k\in\Z.
\end{equation}
Furthermore, if $j$ is odd, then $(x_1,y_1)\in\Z\times\Z$ is a solution to (\ref{dos}) (resp. (\ref{tres})) if and only if (\ref{cuatro})
holds with $k$ odd (resp. even) \cite{S}. Thus
\begin{equation}
\label{cinco}
x+y\d
\end{equation}
is the smallest real number such that $x,y$ are both positive integers and $(x,y)$ is a solution to $(\ref{uno})$. For this reason, we will
refer to $(x,y)$ as the {\em smallest solution} to (\ref{uno}).

For certain values of $d$, a closed formula for (\ref{cinco}) is available. For instance, if $d=e^2+1$, then $(x,y)=(e,1)$ is
the smallest solution to (\ref{dos}), $\d=[e,\overline{2e}]$, and $j=1$. Also, if $d=e^2+2e$, then $(x,y)=(e+1,1)$ is
the smallest solution to (\ref{tres}), $\d=[e,\overline{1,2e}]$, and $j=2$. More generally, if $d=e^2+\frac{2e}{k}$, where $k\in\N$ and $k|2e$, then
\begin{equation}
\label{CF}
(x,y)=(ke+1,k)
\end{equation}
is the smallest solution to (\ref{tres}), $\d=[e,\overline{k,2e}]$, and $j=2$, provided $k\neq 2e$. The following table illustrates a few instances, with $j=2$, of the latter phenomenon.
$$
\begin{array}{c|c|c|c|c}
    e & k & d & x & y \\\hline
    1 & 1 & 3 & 2 & 1 \\\hline
    2 & 1 & 8 & 3 & 1 \\\hline
    2 & 2 & 6 & 5 & 2 \\\hline
    3 & 1 & 15 & 4 & 1 \\\hline
    3 & 2 & 12 & 7 & 2 \\\hline
    3 & 3 & 11 & 10 & 3 \\\hline
    4 & 1 &  24  & 5  & 1 \\\hline
    4 & 2 &  20  &  9 & 2 \\\hline
    4 & 4 &  18  &  17 & 4 \\
  \end{array}
$$

We seek a closed formula, such as (\ref{CF}), rather than algorithm, such as (\ref{zero}), for the smallest solution to (\ref{uno}),
available for a suitable family of values for $d$. The preceding examples suggest that we consider the case

\begin{equation}
\label{siete}
\d=[e,\overline{k,\dots,k,2e}],\text{ with period }j\geq 2.
\end{equation}
Note that $j\geq 2$ implies $k\neq 2e$. However, given $j,k\in\N$, with $j\geq 3$, there may not exist any $d,e\in\N$ such that (\ref{siete}) holds (unlike the case $j=2$). Indeed, let $d\in\N$, not a perfect square, with
$$
\d=[a_0,a_1,a_2,\dots],\text{ with period }j\geq 1,
$$
and set
\begin{equation}
\label{qp0}
\left(
  \begin{array}{cc}
   q_{-2} & p_{-2} \\
    q_{-1} & p_{-1} \\
  \end{array}
\right)=\left(
  \begin{array}{cc}
    1 & 0 \\
    0 & 1 \\
  \end{array}
\right)
\end{equation}
as well as
\begin{equation}
\label{qp1}
\left(
  \begin{array}{cc}
   q_{n-1} & p_{n-1} \\
    q_{n} & p_{n} \\
  \end{array}
\right)=\left(
  \begin{array}{cc}
   0 & 1 \\
    1 & a_n \\
  \end{array}
\right)\cdots \left(
  \begin{array}{cc}
   0 & 1 \\
    1 & a_0 \\
  \end{array}
\right),\quad n\geq 0.
\end{equation}
Then
$$
q_n=a_n q_{n-1}+q_{n-2},\; p_n=a_n p_{n-1}+p_{n-2},\quad n\geq 0,
$$
and, according to (\ref{zero}),
$$
(x,y)=(q_{j-1},p_{j-1})
$$
is the smallest solution to (\ref{uno}). Suppose, if possible, that
$$
\d=[e,\overline{1,1,2e}]
$$
for some $d,e\in\N$. By above,
$$
\begin{array}{c|c|c}
     n & q_n & p_n \\\hline
    -2 & 1 & 0 \\\hline
    -1 & 0 & 1 \\\hline
    0 & 1 & e \\\hline
    1 & 1 & e+1 \\\hline
    2 & 2 & 2e+1 \\\hline
  \end{array}
$$
and $(x,y)=(2e+1,2)$ is a solution to (\ref{uno}). Thus $(2e+1)^2-4d=-1$, which implies
$$
d=\frac{4e^2+4e+2}{4},
$$
against the fact that $d\in\N$. The above discussion leads to the following

\noindent{\bf Goal.} Given $j,k\in\N$, with $j\geq 2$, we wish to find necessary and sufficient conditions for the existence of $d,e\in\N$ such that (\ref{siete}) holds. Moreover, when these
conditions are satisfied, we aim to find all possible values for $d,e$, as well as a closed formula for the smallest solution $(x,y)$
to the Fermat-Pell equation (\ref{uno}).

For instance, when $k\in\N$, we look for all $e,d\in\N$ so that
\begin{equation}
\label{ocho}
\d=[e,\overline{k,k,2e}],\text{ with period }j=3,
\end{equation}
as well as the smallest solution to (\ref{uno}). Assuming (\ref{ocho}) holds, the algorithm (\ref{zero}) produces
$$
\begin{array}{c|c|c}
     n & q_n & p_n \\\hline
    -2 & 1 & 0 \\\hline
    -1 & 0 & 1 \\\hline
    0 & 1 & e \\\hline
    1 & k & ke+1 \\\hline
    2 & k^2+1 & (k^2+1)2e+k \\\hline
  \end{array}
$$
where $q_n,p_n$ are defined in (\ref{qp0}) and (\ref{qp1}), so
\begin{equation}
\label{nueve}
(x,y)=((k^2+1)e+k,k^2+1)
\end{equation}
is the smallest solution to (\ref{uno}). Thus
$$
((k^2+1)e+k)^2-d(k^2+1)^2=-1,
$$
which implies
\begin{equation}
\label{diez}
d=\frac{((k^2+1)e+k)^2+1}{(k^2+1)^2}=e^2+\frac{2(k^2+1)ke+k^2+1}{(k^2+1)^2}=e^2+\frac{2ke+1}{k^2+1}.
\end{equation}
Now
$$
2ke\equiv -1\mod k^2+1
$$
is solvable if and only if $k$ is even, in which case there is a unique solution $e\equiv \frac{k}{2} \mod k^2+1$, whence
\begin{equation}
\label{once}
e=\frac{k}{2}+\ell (k^2+1),\quad \ell\geq 1\quad (\text{to ensure period }j=3).
\end{equation}
Suppose, conversely, that $x,y,d,e$ are as indicated in (\ref{nueve})-(\ref{once}).
Then, by construction, we have $x/y=[e,k,k]$. Since $x^2-dy^2=-1$, it follows from \cite[Exercises 7.7.15, 7.7.17]{S} that
$\d=[e,\overline{k,k,2e}]$, with period $j=3$ as $2e\neq k$. (There is a typing error at the beginning of \cite[Exercise 7.7.17]{S}: It should
say $X^2-dY^2=\pm 1$ instead of $X^2-dY^2=1$.)

The following table lists a few examples.
$$
\begin{array}{c|c|c|c|c}
    e & k & d & x & y \\\hline
    6 & 2 & 41 & 32 & 5 \\\hline
    11 & 2 & 130 & 57 & 5 \\\hline
    16 & 2 & 269 & 82 & 5 \\\hline
    19 & 4 & 370 & 327 & 17 \\\hline
    36 & 4 & 1313 & 616 & 17 \\\hline
    40 & 6 & 1613 & 1486 & 37 \\\hline
    69 & 8 &  4778  & 4493  & 65 \\
    \end{array}
$$

\section{The general case}

Given $k\in\N$ and $j=m+1\geq 2$, suppose that (\ref{siete}) holds, and let $(x,y)$ be the smallest solution to (\ref{uno}).
Set
\begin{equation}
\label{def}
f_{-2}=1,\; f_{-1}=0,\; f_n=kf_{n-1}+f_{n-2},\quad n\geq 0.
\end{equation}
Then the algorithm (\ref{zero}) yields
\begin{equation}
\label{cuad}
(x,y)=(f_m+ef_{m-1},f_m).
\end{equation}

We digress here to record a few properties of the sequence (\ref{def}). We readily see by induction that
\begin{equation}
\label{gcd}
\gcd(f_n,f_{n-1})=1,\quad n\geq -1
\end{equation}
as well as
$$
\left(
  \begin{array}{cc}
    0 & 1 \\
    1 & k \\
  \end{array}
\right)^n=\left(
  \begin{array}{cc}
    f_{n-2} & f_{n-1} \\
    f_{n-1} & f_{n} \\
  \end{array}
\right),\quad n\geq 1,
$$
which implies
\begin{equation}
\label{cassini}
f_{n-1}^2+(-1)^n=f_n f_{n-2},\quad n\geq 1 \text{ (Cassini identity).}
\end{equation}
Arguing by induction, and making use of (\ref{def}) and (\ref{cassini}), we easily derive
\begin{equation}
\label{cassini2}
(-1)^n f_{n-2}f_{n-1}\equiv k\mod f_n,\quad n\geq 2.
\end{equation}

The following table displays $f_n$ for all $-2\leq n\leq 11$.
$$
\begin{array}{c|c}
   n & f_n\\\hline
  -2 & 1 \\\hline
  -1 & 0 \\\hline
  0 & 1 \\\hline
  1 & k \\\hline
  2 & k^2+1 \\\hline
  3 & k^3+2k \\\hline
  4 & k^4+3k^2+1 \\\hline
  5 & k^5+4k^3+3k \\\hline
  6 & k^6+5k^4+6k^2+1 \\\hline
  7 & k^7+6k^5+10k^3+4k \\\hline
  8 & k^8+7k^6+15k^4+10k^2+1 \\\hline
  9 & k^9+8k^7+21k^5+20k^3+5k \\\hline
  10 & k^{10}+9k^8+28k^6+35k^4+15k^2+1 \\\hline
  11 & k^{11}+10k^9+36k^7+56k^5+35k^3+6k  \\
\end{array}
$$
As the table suggests, it is easily shown by induction that
$$
f_n=a_nk^n+a_{n-2}k^{n-2}+a_{n-4}k^{n-4}+\cdots,\quad n\geq 0,
$$
where
$$
a_{n-2i}={{1}\choose{0}}+{{i}\choose{1}}+{{i+1}\choose{2}}+\cdots+{{n-i-1}\choose{n-2i}}={{1}\choose{0}}+\underset{1\leq s\leq n-2i}\sum
{{s+(i-1)}\choose{s}}.
$$

We now resume our prior discussion. Since (\ref{cuad}) is a solution to (\ref{uno}), we have
$$
d=\frac{(f_m e+f_{m-1})^2+(-1)^m}{f_m^2}=e^2+\frac{2f_mf_{m-1}e+f_{m-1}^2+(-1)^m}{f_m^2}.
$$
By virtue of (\ref{cassini}), the above may rewritten as follows:
\begin{equation}
\label{u}
d=e^2+\frac{2f_{m-1}e+f_{m-2}}{f_m}.
\end{equation}
The congruence equation
\begin{equation}
\label{congo}
2f_{m-1}e\equiv -f_{m-2}\mod f_m
\end{equation}
is solvable if and only if
$$
\gcd(f_m,2f_{m-1})|f_{m-2}.
$$
By (\ref{gcd}), the above translates into
$$
\gcd(f_m,2)|f_{m-2},
$$
which is true except only when $f_m$ is even and $f_{m-2}$ is odd. Suppose first $k$ is even. Then $f_m,f_{m-2}$
are both even if $m$ is odd (no constant term) and both odd if $m$ is even (constant term 1). Suppose next $k$ is
odd. Then (\ref{def}) yields the following parity for $f_n$:
$$
\begin{array}{ccccccc}
  f_0 & f_1 & f_2 & f_3 & f_4 & f_5 & \dots \\
  \text{odd} & \text{odd} & \text{even} & \text{odd} & \text{odd} & \text{even} & \dots
\end{array}
$$
We deduce that
$$
f_m\text{ is even and }f_{m-2}\text{ is odd }\Leftrightarrow k\text{ is odd and }m\equiv 2\mod 3.
$$
Suppose henceforth that $k$ is even or $m\not\equiv 2\mod 3$. Since $\gcd(f_m,f_{m+1})=1$, we may multiply (\ref{congo}) by $f_{m+1}$
and obtain the equivalent equation
$$
2f_{m-1}f_{m+1}e\equiv -f_{m-2}f_{m+1}\mod f_m.
$$
By (\ref{cassini}), this translates into
$$
2(f_m^2+(-1)^{m+1})e\equiv -f_{m-2}f_{m+1}\mod f_m,
$$
that is,
$$
2e\equiv (-1)^m f_{m-2}f_{m+1}\equiv (-1)^m f_{m-2}(kf_m+f_{m-1})\equiv (-1)^m f_{m-2}f_{m-1}\mod f_m.
$$
In light of (\ref{cassini2}), this is equivalent to
\begin{equation}
\label{magen}
2e\equiv k\mod f_m.
\end{equation}
Three cases arise.

\noindent{\sc Case 1.} $k$ is even and $f_m$ is odd. Then all solutions to (\ref{magen}) are of the form
\begin{equation}
\label{v}
e=\frac{k}{2}+\ell f_m,\quad \ell\geq 1\quad (\text{to ensure period }j).
\end{equation}
By (\ref{u}) and (\ref{v}), we see that 
$$
d=e^2+\frac{(k+2\ell f_m)f_{m-1}+f_{m-2}}{f_m}=e^2+\frac{2\ell f_mf_{m-1}+kf_{m-1}+f_{m-2}}{f_m}=
e^2+\frac{2\ell f_mf_{m-1}+f_m}{f_m},
$$
so
\begin{equation}
\label{d1}
d=e^2+2\ell f_{m-1}+1,\quad \ell\geq 1.
\end{equation}
\noindent{\sc Case 2.} $k$ is even and $f_m$ is even. Then all solutions to (\ref{magen}) are of the form
\begin{equation}
\label{v2}
e=\frac{k}{2}+\ell \frac{f_m}{2},\quad \ell\geq 1,
\end{equation}
in which case the above calculations yield
\begin{equation}
\label{d2}
d=e^2+\ell f_{m-1}+1,\quad \ell\geq 1.
\end{equation}
\noindent{\sc Case 3.} $k$ is odd. In this case $m\not\equiv 2\mod 3$, so $f_m$ is also odd. Then $k+f_m$ is even, and
all solutions to (\ref{magen}) are of the form
\begin{equation}
\label{v3}
e=\frac{k+f_m}{2}+\ell f_m,\quad \ell\geq 0,
\end{equation}
in which case the above calculations yield
\begin{equation}
\label{d3}
d=e^2+(2\ell+1)f_{m-1}+1,\quad \ell\geq 0.
\end{equation}

Suppose, conversely, that $x,y$ are as in (\ref{cuad}) and that $d,e$ are as indicated in (\ref{v})-(\ref{d3}).
Then, by construction, we have $x/y=[e,k,\dots,k]$, where $k$ is repeated $m$ times. Since $x^2-dy^2=-1$, it follows from \cite[Exercises 7.7.15, 7.7.17]{S} that $\d=[e,\overline{k,\dots, k, 2e}]$, with period $j$ as $2e\neq k$.

We have proven the following

\noindent{\bf Theorem. }{\em Let $j,k\in\N$, with $j=m+1\geq 2$, and let $f_{-2},f_{-1},f_0,\dots$ be defined as in (\ref{def}).
Then there exist $d,e\in\N$ such that
$$
\d=[e,\overline{k,\dots,k,2e}],\text{ with period }j\geq 2,
$$
if and only if $k$ is even or $m\not\equiv 2\mod 3$, in which case all such $d,e$ are given by (\ref{v}) and (\ref{d1}) if $k$ is even
and $f_m$ is odd, by (\ref{v2}) and (\ref{d2}) if $k$ is even
and $f_m$ is even, and by (\ref{v3}) and (\ref{d3}) if $k$ is odd; moreover, the smallest solution to (\ref{uno}) is as indicated in (\ref{cuad}).}

The following table illustrates a few instances of the above result.

$$\begin{array}{c|c|c|c|c|c|c|c}
    m & k & f_{m-1} & f_m & \ell & e & d & f_me+f_{m-1} \\\hline
    3 & 1 & 2 & 3      &  0  & 2   & 7      & 8 \\\hline
    3 & 1 & 2 & 3      &  1  & 5   & 32     & 17 \\\hline
    3 & 1 & 2 & 3      &  2  & 8   & 75     & 26\\\hline
    3 & 1 & 2 & 3      &  3  & 11  & 136    & 35 \\\hline
    3 & 2 & 5 & 12     &  4  & 25  & 646    & 305 \\\hline
    3 & 3 & 10 & 33    &  0  & 18  & 335    & 604  \\\hline
    3 & 4 & 17 & 72    &  1  & 38  & 1462   & 2753 \\\hline
    3 & 5 & 26 & 135   &  0  & 70  & 4927   & 9476 \\\hline
    4 & 1 & 3 & 5      &  0  & 3   & 13     & 18 \\\hline
    4 & 1 & 3 & 5      &  1  & 8   & 74     & 43 \\\hline
    4 & 2 & 12 & 29    &  1  & 30  & 925    & 882 \\\hline
    4 & 3 & 33 & 109   &  0  & 56  & 3170   & 6137 \\\hline
    5 & 2 & 29 & 70    &  2  & 71  & 5100   & 4999 \\\hline
    5 & 2 & 29 & 70    &  4  & 141 & 19998  & 9899 \\\hline
    5 & 2 & 29 & 70    &  6  & 211 & 44696  & 14799 \\\hline
    5 & 4 & 305 & 1292 &  2  &1294 & 1675047& 1672153 \\\hline
    6 & 1 & 8 & 13     &  0  &  7  & 58     & 99 \\\hline
    6 & 1 & 8 & 13     &  1  & 20  & 425    & 268 \\\hline
    6 & 2 & 70 & 169   &  1  & 170 & 29041  & 28800 \\\hline
    7 &  1 & 13 & 21   &  0  & 11  & 135    & 244 \\\hline
    7 &  2 & 169 & 408 &  2  & 409 & 167620 & 167041 \\\hline
    9 &  1 & 34 & 55   &  0  & 28  & 819    & 1574 \\\hline
    10 & 1 & 55 & 89   &  0  & 45  & 2081   & 4060 \\\hline
    12 & 1 & 144 & 233 &  0  & 117 & 13834  & 27405 \\\hline
    13 & 1 & 233 & 377 &  0  & 189 & 35955  & 71486 \\
\end{array}
$$

\end{document}